\documentclass[10pt]{amsart}
\usepackage{amscd}

\newtheorem{thm}{Theorem}[section]
\theoremstyle{definition}
\newtheorem{defn}{Definition}[section]
\theoremstyle{remark}

\begin{document}
\title{Subspaces of Knot spaces}
\author{Craig Benham}
\author{Xiao-Song Lin}
\author{David Miller}

\begin{abstract}
The inclusion of the space of all knots of a prescribed writhe in a
particular isotopy class into the space of all knots in that isotopy class
is a weak homotopy equivalence.
\end{abstract}

\maketitle
\section{Introduction}

In \cite{benhammiller}, 
it is shown that the set of knots in any fixed isotopy class that
have a prescribed writhe is path-connected. This result can be interpreted
in the following manner. Let 
\begin{equation*}
i^{\mathcal{S}}:\,\,\mathcal{S}_{w}{\hookrightarrow }\mathcal{M}_{c}
\end{equation*}
represent the inclusion of the space of a prescribed writhe $w$ into the
space of all knots in some isotopy class. Then the results of 
\cite{benhammiller} can be
understood as stating that this map induces an isomorphism on $\pi _{0}$. In
this paper, we prove that this inclusion is, in fact, a weak homotopy
equivalence, inducing isomorphisms on $\pi _{i}$ for all $i$.

The fact that $(i^{\mathcal{S}})_{*}:\pi _{n}(\mathcal{S}_{w})\rightarrow
\pi _{n}(\mathcal{M}_{c})$ is an isomorphism will be proved using an
extension of the methods of \cite{benhammiller}.

\section{Notation, definitions and statement of theorems}

Let $S^{1}=[0,1]\big/0\sim 1$. We define a closed, space curve to be a
smooth ($C^{1}$) embedding 
\begin{equation*}
\gamma :\,\,S^{1}\hookrightarrow \mathbb R^{3}.
\end{equation*}
In general, we stipulate that $\gamma ^{\prime }(s)\neq 0$. The image of the
map $\gamma $ is referred to as a smooth knot.

\begin{defn}
The writhe of the embedding $\gamma $ is given by 
\begin{equation*}
\mathbf{Wr}(\gamma )=\iint_{\gamma (S^{1})\times \gamma (S^{1})-\Delta }dA
\end{equation*}
where $dA$ is the pull back of the area element on the unit sphere under the
map which assigns to every pair $(\gamma (s),\gamma (s^{\prime }))\in \gamma
(S^{1})\times \gamma (S^{1})-\Delta $ the unit vector from $\gamma (s)$ to $%
\gamma (s^{\prime })$.
\end{defn}

By definition, the writhe is a geometric property of the image of $\gamma $,
a property of the knot. in particular, it is independent of
parameterization. We can thus refer alternatively to the writhe, $\mathbf{Wr}%
(K)$, of the knot $K$. We let $\mathcal{M}_{c}$ denote an arbitrary
component of the space of knots, and $\mathcal{S}_{\omega }$ denote the
subset of $\mathcal{M}_{c}$ of curves with a particular writhe, $\omega $.

\begin{thm}
The induced map 
\begin{equation*}
(i^{\mathcal{S}})_{*}:\,\,\pi _{n}(\mathcal{S}_{w},\gamma _{0})\rightarrow \pi
_{n}(\mathcal{M}_{c},\gamma _{0})
\end{equation*}
is an isomorphism for all $n$,where the basepoint $\gamma _{0}$ is a smooth
knot of writhe $\omega $.
\end{thm}

The proof of theorem 2.1 will extend the results and proofs of 
\cite{benhammiller}, which
can be understood as theorem 2.1 in the case $n=0$. In particular, we need
the following theorem of Fuller.

\begin{thm}[Fuller \cite{fuller1}]
The unit tangents $T(s)$ to a smooth, closed space curve,$\gamma $, trace
out, if their starting points are translated to the origin, a closed curve
on the unit sphere. Let $A$ be the area on the unit sphere enclosed by this
curve. Then 
\begin{equation*}
1+\mathbf{Wr}(\gamma )=\frac{1}{2\pi }A,\newline
\text{\textrm{modulo}}\ 2.
\end{equation*}
\end{thm}

\section{Construction of fixed-writhe families of embeddings}

In this section, we will generalize the results of \cite{benhammiller}
 in two different
settings, in order to analyze the higher connectivity of the sets of fixed
writhe within components of the space of knots.

First consider any family of embedded curves 
\begin{equation*}
\gamma _x:\,\,S^1\hookrightarrow\mathbb R^3,\qquad\qquad x\in(S^n,x_0)
\end{equation*}
indexed by the based $n$-dimensional sphere such that the map 
\begin{equation*}
\Gamma :\,\,S^1\times S^n\rightarrow\mathbb R^3\times S^n \qquad\qquad
\Gamma(s,x)\overset{\mathbf{def.}}{=}( \gamma_x(s),x)
\end{equation*}
is continuous. From $\Gamma $ we can construct 
\begin{equation*}
\Psi :\,\,S^n\rightarrow \mathcal{M}_c,\qquad\qquad\Psi (x) \overset{\mathbf{def.%
}}{=}\gamma_x
\end{equation*}

Let 
\begin{equation*}
\omega=\mathbf{Wr}(\gamma_{x_0})
\end{equation*}
Our goal is to construct a family of embeddings $\bar\gamma_x$ such that $%
\mathbf{Wr}(\gamma_{x})=\omega$ for all $x\in S^n$. To accomplish this goal,
we must first alter each embedding in the original family in order to allow
for a correction in the writhe.

We choose once and for all, a point $s_0\in S^1$, $0\neq s_0\neq 1$. Let $%
T_x(s_0)$ be the unit tangent to $\gamma_x$ at $s_0$. We reparameterize the
maps $\gamma_x$ such that for each $x\in S^n$, there is an interval $%
[s_1,s_2]$ containing $s_0$ and for any $s\in [s_1,s_2]$, $%
\gamma_x(s)=\gamma_x(s_0)$. This can be accomplished by reparameterizing $%
\Gamma$ so that it is constant on $[s_1,s_2]\times S^n\subset S^1\times S^n$
and not on any larger closed set, i.e., $\Gamma^{-1}(\gamma_x(s_0))=
[s_1,s_2]\times I$. Then for each $x\in [0,1]$, let $B_{x,\epsilon}\subset%
\mathbb R^3$ be a closed ball of radius $\epsilon$ with $\gamma_x(s_0)$ at
its center. $\epsilon$ can be chosen arbitrarily small, but is fixed for all 
$x$. Let $\phi_x:B_{x,\epsilon}\rightarrow B_\epsilon\hookrightarrow\mathbb %
R^3$ be a linear, co-ordinate chart for $B_{x,\epsilon}$, i.e., a
translation map taking the point $\gamma_x(s_0)$ to the origin, and the ball 
$B_{x,\epsilon}$ to a ball $B_\epsilon$ of radius $\epsilon$ centered at the
origin.

We construct for each $x$ a map $\sigma _{x}:(B_{\epsilon }-\{\mathbf{0}%
\})\longrightarrow (B_{\epsilon }-\{\mathbf{0}\})$ as follows: 
\begin{equation*}
\sigma _{x}(z)=\cases\big[(1-\frac{3|x|}{2})|z|+\frac{3|x|\epsilon }{2}\big]%
\frac{z}{|z|},\qquad 0\leq |x|\leq \frac{1}{3} \\[2mm]
\big[\frac{1}{2}|z|+\frac{\epsilon }{2}\big]\frac{z}{|z|},\qquad \qquad
\qquad |x|\geq \frac{1}{3}.\endcases
\end{equation*}
These linear maps have a number of important properties. First note that $%
\sigma _{x_{0}}=\mathbf{id}_{B_{\epsilon }}$. Note also that for each $x$, $%
\sigma _{x}|_{S_{\epsilon }}=\mathbf{id}_{S_{\epsilon }}$. Most importantly,
these maps have the property that they ``push'' the contents of $B_{\epsilon
}-\{\mathbf{0}\}$ towards the boundary in a linear fashion, creating an
``empty space'' around the origin. When $x$ is close to $x_{0}$, the
contents are not pushed as far.

We next assemble a map $\tilde\sigma :\mathbb R^3\rightarrow\mathbb R^3$ by 
\begin{equation*}
\tilde\sigma_x(z)=\cases z,\qquad\qquad\qquad z\in\mathbb %
R^3-(B_{x,\epsilon}- \{\gamma_x(s_0)\}) \\[1.3mm]
\big(\phi_x^{-1}\circ\sigma_x\circ \phi_x \big)(z),\qquad z\in
B_{x,\epsilon}-\{\gamma_x(s_0)\}. \endcases
\end{equation*}
If, for $x\neq x_0$, we think of our knot, $\mathbf{im}(\gamma_x)$, lying in 
$\mathbb R^3$, its image under $\tilde\sigma_x$ would no longer be a closed
curve. There would be a gap between the two points where the curve touched a
certain small sphere around $\gamma_x(s_0)$. To rectify this problem, we
construct an altered isotopy, 
\begin{equation*}
\tilde \gamma_x(s)= \cases (\tilde\sigma_x\circ \gamma_x)(s),
\qquad\qquad\qquad s\in (S^1-[s_1,s_2]) \\[1.3mm]
\phi_x^{-1} \big(\frac{2\epsilon}{s_2-s_1} [s-\frac{s_1+s_2}{2}] (T_x(s_0))%
\big), \qquad s\in [s_1, s_2]. \endcases
\end{equation*}
The definition of $\tilde\gamma_x$ on $[s_1,s_2]$ has the effect of
inserting a straight line segment in the gap left by $\tilde\sigma_x$, and
connecting smoothly with $\tilde\gamma_x(s_1)$ and $\tilde\gamma_x(s_2)$ at
the ends.

The final step in constructing a family of embeddings of fixed writhe is to
replace the altered portion of the curves with helical segments which will
provide the appropriate correction in the writhe. By Theorem 2.2, the writhe
of a knot is determined modulo an integral indeterminacy by the area on the
unit sphere enclosed by the tangent indicatrix, the curve that the unit
tangents to the knot trace out. We will use this fact to construct a local
correction to the altered isotopy above which will precisely fix the writhe
throughout the isotopy. Let $M=\underset{x}{\text{max}}\big(\mathbf{Wr}(%
\tilde{\gamma}_{x})\big)$ and $m=\underset{x}{\text{min}}\big(\mathbf{Wr}(%
\tilde{\gamma}_{x})\big)$. Let $n$ be the smallest integer strictly greater
than $M-m$. Let $w(x)=\omega -\mathbf{Wr}(\tilde{\gamma}_{x})$. Note that $%
|w(x)|<n$. We shall modify each embedding $\tilde{\gamma}_{x}:S^{1}%
\hookrightarrow \mathbb R^{3}$ on the interval $[s_{1},s_{2}]$. Let 
\begin{equation*}
s_{3}=s_{1}+\frac{s_{2}-s_{1}}{4},\qquad \qquad \qquad s_{4}=s_{2}-\frac{%
s_{2}-s_{1}}{4}
\end{equation*}
On $[s_{3},s_{4}]$, and for each $x\in S^{n}$ we define the following helix
function 
\begin{align*}
\tau _{x}(s)& =\big(r(x)\cos (C(x)s),r(x)\sin (C(x)s), \\[2mm]
& p(x)\frac{n}{s_{4}-s_{3}}(s-\frac{s_{4}+s_{3}}{2})\big),\qquad \text{%
\textrm{where}} \\
C(x)& =-\text{\textrm{sgn}}(w(x))\frac{2\pi n}{s_{4}-s_{3}} \\[2mm]
r(x)& =S(x)\frac{\epsilon }{4\pi n}\sqrt{\frac{2|w(x)|}{n}-\frac{w(x)^{2}}{%
n^{2}}} \\[2mm]
p(x)& =S(x)\frac{\epsilon }{2n}(1-\frac{|w(x)|}{n}),
\end{align*}
and where $S(x)$ is a scale function given by 
\begin{equation*}
S(x)=\cases {3x},\qquad 0\leq |x|\leq \frac{1}{3} \\[1.5mm]
1,\qquad |x|>\frac{1}{3}\endcases
\end{equation*}
Thus the image of $\tau _{x}$ is a helix with precisely $n$ turns. The pitch
and radius of the helix are determined by $w(x)$, but each turn of the helix
is contained in a ball of radius $S(x)\frac{\epsilon }{2n}$. If $w(x)$ goes
to zero, the helix degenerates into a vertical line of length $\frac{%
\epsilon }{2}$, from $(0,0,-\frac{\epsilon }{4})$ to $(0,0,\frac{\epsilon }{4%
})$. If $|w(x)|$ gets close to $n$, the pitch angle becomes small.

We must now connect the image of $\tau _{x}$ with the south and north poles
of the sphere of radius $S(x)\frac{\epsilon }{2}$ centered at the origin. To
do this we consider the plane spanned by the unit vectors 
\begin{equation*}
\mathbf{v}_{1}=(0,0,1)\qquad \text{\textrm{and}}\qquad \mathbf{v}_{2}=\frac{%
\tau _{x}^{\prime }(s_{3})}{|\tau _{x}^{\prime }(s_{3})|}=\frac{\tau
_{x}^{\prime }(s_{4})}{|\tau _{x}^{\prime }(s_{4})|}.
\end{equation*}
We then choose $\iota _{x}:[s_{1},s_{3}]\hookrightarrow (B_{S(x)\frac{%
\epsilon }{2}}^{-}-B_{\Vert \tau _{x}(s_{3})\Vert }^{-})$ be any smooth
curve in the $\mathbf{v}_{1}$-$\mathbf{v}_{2}$ plane such that 
\begin{eqnarray*}
\iota _{x}(s_{1}) &=&(0,0,-\frac{\epsilon }{2}),\qquad \frac{\iota
_{x}^{\prime }(s_{1})}{|\iota _{x}^{\prime }(s_{1})|}=(0,0,1),\qquad \\
\iota _{x}(s_{3}) &=&\tau _{x}(s_{3}),\qquad \frac{\iota _{x}^{\prime
}(s_{3})}{|\iota _{x}^{\prime }(s_{3})|}=\frac{\tau _{x}^{\prime }(s_{3})}{%
|\tau _{x}^{\prime }(s_{3})|}
\end{eqnarray*}
There are clearly many such curves, as we are merely requiring that $\iota $
be a curve on a plane with specified end-points, and specified tangent
directions at those end points. (In the special case where $\mathbf{v}_{1}=%
\mathbf{v}_{2}$, the situation is even simpler. We merely connect the two
end-points with a vertical line.) Our next step is to connect $\tau _{x}$ to
the north pole by defining $\xi _{x}:[s_{3},s_{4}]\hookrightarrow (B_{S(x)%
\frac{\epsilon }{2}}^{+}-B_{\Vert \tau _{x}(s_{4})\Vert }^{+})$ by 
\begin{equation*}
\xi _{x}(s)=r\big(\iota _{x}(s_{3}-(s-s_{4}))\big)
\end{equation*}
where $r$ is reflection about the $x_{1}x_{2}$-plane. We now assemble a map 
\begin{align*}
& \eta _{x}:\,\,[s_{1},s_{2}]\hookrightarrow B_{S(x)\frac{\epsilon }{2}}\ \ 
\text{\textrm{as follows:}} \\[1.3mm]
& \eta _{x}(s)=\cases\iota _{x}(s),\qquad s_{1}\leq s<s_{3}, \\
\tau _{x}(s),\qquad s_{3}& \leq s\leq s_{4}, \\
\xi _{x}(s),\qquad s_{4}& <s<s_{2}, \\
& \endcases
\end{align*}
Let $\theta \in SO(3)$ be any rotation map which takes the unit vector $%
(0,0,1)$ to the unit vector $T_{x}(s_{0})$. Finally, we assemble the map 
\begin{equation*}
\bar{\gamma}_{x}(s)=\cases\tilde{\gamma}_{x}(s),\qquad s\in
(S^{1}-[s_{1},s_{2}]) \\[1.3mm]
\phi _{x}^{-1}\big(\theta (\eta _{x}(s))\big),\qquad s\in [s_{1},s_{2}]%
\endcases
\end{equation*}

We shall prove that $\mathbf{Wr}(\bar{\gamma}_{x})=\omega $ for all $x\in
S^{n}$. We use theorem 2.2, comparing the area that the tangent indicatrix
to $\bar{\gamma}_{x}$ traces out to that traced out by $\tilde{\gamma}_{x}$.
The difference will be equal to $w(x)$. We note first that $\gamma _{x}$ was
parameterized to be constant on the interval $[s_{1},s_{2}]$, so that the
tangent indicatrix curve is not well-defined. However, our analysis was
conducted for $\tilde{\gamma}_{x}$. This curve behaved linearly on $%
[s_{1},s_{2}]$. Thus the tangent indicatrix curve for $\tilde{\gamma}_{x}$
on the unit sphere, $\widetilde{T}_{x}(s)$, is constant on $[s_{1},s_{2}]$.
In defining $\bar{\gamma}$ we changed the behavior of the curve on the
interval $[s_{1},s_{2}]$ in a precise manner. The effect this has on the
tangent indicatrix and the area it encloses can be easily determined. The
intervals $[s_{1},s_{3}]$, and $[s_{4},s_{2}]$, determined respectively by $%
\iota _{x}$ and $\xi _{x}$ have closely related tangent behavior. As $s$
goes from $s_{1}$ to $s_{3}$, the tangent indicatrix will travel along some
curve from $\bar{T}_{x}(s_{1})$ to $\bar{T}_{x}(s_{3})$. While $s$ goes from 
$s_{4}$ to $s_{2}$, the indicatrix of $\bar{T}_{x}(s)$ will travel back
along the same curve. Thus, these intervals together do not contribute to
the area enclosed by the tangent indicatrix. It only remains to examine the
tangent indicatrix curve on the interval $[s_{3},s_{4}]$. On this interval, $%
\bar{\gamma}_{x}$ is a helix which winds $n$ times, and by virtue of the
definition of $\tau _{x}$ has pitch angle $\psi _{x}$ with the property that 
\begin{align*}
& \sin \psi _{x}=1-\frac{w(x)}{n}\ \ \ \text{\textrm{which means that}} \\
& w(x)=n(1-\sin \psi _{x})
\end{align*}
But this last expression is precisely $n$ times the signed area of the \emph{%
cap} on the unit sphere that the unit tangents to an $n$-turn helix of pitch
angle $\psi _{x}$ trace out $n$ times. Thus the insertion of the helix $\tau
_{x}$ contributes precisely $w(x)$ to the area traced out by the tangent
indicatrix. (Since the writhe changes continuously through the isotopy, the
indeterminacy in Theorem 2.2 does not arise.) Thus, we see that 
\begin{equation*}
\mathbf{Wr}(\bar{\gamma}_{x})=\mathbf{Wr}(\tilde{\gamma}_{x})+w(x){=}\omega .
\end{equation*}

In the next setting to which we will generalize the results of 
\cite{benhammiller}, we
consider families of embedded curves 
\begin{equation*}
\lambda_x:\,\,S^1\hookrightarrow\mathbb R^3,\qquad\qquad x\in (S^n\times I,x_0)
\end{equation*}
such that 
\begin{equation*}
\Lambda:\,\,S^1\times (S^n\times I)\rightarrow \mathbb R^3\times(S^n\times
I),\qquad\qquad \Lambda(s,x)\overset{\mathbf{def.}}{=}(\lambda_x(s),x).
\end{equation*}
In addition, we stipulate that if $x\in(S^n\times\{0\})$ or $%
x\in(S^n\times\{1\})$, that $\mathbf{Wr}(\lambda_x)=\omega$. Our goal is,
again, to construct a family of embeddings $\bar\lambda_x$ such that $%
\mathbf{Wr}(\bar\lambda_x)=\omega$, for all $x\in (S^n\times I)$. From $%
\Lambda$, we can construct 
\begin{equation*}
\Upsilon :\,\,S^n\times I\rightarrow\mathcal{M}_c.\qquad\qquad\Upsilon(x) 
\overset{\mathbf{def.}}{=}\lambda_x.
\end{equation*}
The existence of the curves $\bar\lambda_x$ is proved with exactly the
construction of \cite{benhammiller} crossed with $S^n$.

\section{Proof of theorem 2.1}

We first show that the map $(i^{\mathcal{S}})_{*}$ is onto. We first define 
\begin{equation*}
\overline{\Gamma }:\,\,S^{1}\times S^{n}\rightarrow \mathbb R^{3}\times
S^{n}\qquad \qquad \overline{\Gamma }(s,x)\overset{\mathbf{def.}}{=}(\bar{%
\gamma}_{x}(s),x)
\end{equation*}
and 
\begin{equation*}
\overline{\Psi }:\,\,S^{n}\rightarrow \mathcal{M}_{c},\qquad \qquad \overline{%
\Psi }(x)\overset{\mathbf{def.}}{=}\bar{\gamma}_{x}
\end{equation*}
We define $\tilde{\Gamma}$ and $\tilde{\Psi}$ similarly. We then define a
homotopy 
\begin{equation*}
\Omega :\,\,S^{n}\times I\rightarrow \mathcal{M}_{c}\times I
\end{equation*}
as follows. Let $\Omega _{0}\overset{\mathbf{def.}}{=}\Psi $ and $\Omega _{1}%
\overset{\mathbf{def.}}{=}\overline{\Psi }$. For $\frac{1}{2}\leq t\leq 1$
we define $\Omega _{t}$ to be $\overline{\Psi }$, except that in the
definition of $\bar{\gamma}_{x}$, we replace $r(x)$ with $(2t-1)r(x)$. This
has the effect of shrinking the width of the inserted helix as $t$ decreases
from $1$ to $\frac{1}{2}$, so that $\Omega _{\frac{1}{2}}=\tilde{\Psi}$. For 
$0\leq t<\frac{1}{2}$, we define $\Omega _{t}$ to be $\tilde{\Psi}$, except
that in the definition of $\tilde{\gamma}_{x}$, we replace $\epsilon $ with $%
2t\epsilon $. We then see that $\Omega $ is a homotopy between $\Psi $ and $%
\overline{\Psi }$.

If we take any element $y\in \pi _{n}(\mathcal{M}_{c})$ and write $y=[\Psi ]$%
, the existence of $\Omega $ implies that $y=[\overline{\Psi }]$. But since 
\begin{equation*}
\overline{\Psi }:\,\,S^{n}\rightarrow \mathcal{S}_{w},
\end{equation*}
we see that $y\in i_{*}^{\mathcal{S}}(\pi _{n}(\mathcal{S}_{w}))$. Since the
choice of $y$ was arbitrary, $i_{*}^{\mathcal{S}}$ is onto.

We now show that the map $(i^{\mathcal{S}})_{*}$ is one-one. We define 
\begin{equation*}
\overline{\Lambda }:\,\,S^{1}\times (S^{n}\times I)\rightarrow \ \mathbb %
R^{3}\times (S^{n}\times I),\qquad \qquad \overline{\Lambda }(s,x)\overset{%
\mathbf{def.}}{=}(\bar{\lambda}_{x}(s),x).
\end{equation*}
and we define 
\begin{equation*}
\overline{\Upsilon }:\,\,(S^{n}\times I)\rightarrow \mathcal{M}_{c},\qquad
\qquad \overline{\Upsilon }(x)\overset{\mathbf{def.}}{=}\bar{\lambda}_{x}
\end{equation*}
We define $\tilde{\Lambda}$ and $\tilde{\Upsilon}$ similarly. We define a
homotopy 
\begin{equation*}
\Phi :\,\,(S^{n}\times I)\times I\rightarrow \mathcal{M}_{c}\times I,
\end{equation*}
as follows. For $\frac{1}{2}\leq t\leq 1$, let $\Phi _{t}$ be the same as $%
\overline{\Upsilon }$ except that in the definition of $\bar{\lambda}_{x}$,
we replace $r(x)$ with $(2t-1)r(x)$. Thus, $\Phi _{1}=\overline{\Upsilon }$
and $\Phi _{\frac{1}{2}}=\tilde{\Upsilon}$. For $0\leq t<\frac{1}{2}$, we
define $\Phi _{t}$ to be $\tilde{\Upsilon}$, except that in the definition
of $\tilde{\lambda}_{x}$, we replace $\epsilon $ by $2t\epsilon $. We then
see that $\Phi $ is a homotopy between $\Upsilon $ and $\overline{\Upsilon
}$

If we take two elements $y,z\in \pi _{n}(\mathcal{S}_{w})$, represented
respectively by $\Theta $ and $\Sigma $, such that 
\begin{equation*}
(i^{\mathcal{S}})_{*}(y)=(i^{\mathcal{S}})_{*}(z),
\end{equation*}
there exists a homotopy which we will denote by $\Lambda $, between $\Theta $
and $\Sigma $. But since $\Lambda $ can be homotoped to $\overline{\Lambda }$
via $\Phi $ above, we see that $\overline{\Lambda }$ is a homotopy between $%
\Theta $ and $\Sigma $ within $\mathcal{S}_{w}$. Thus $y=z$ and $(i^{%
\mathcal{S}})_{*}$ is one-one.

\section{Remarks}

It should be noted that the Smale conjecture, proved in \cite{hatcher1}, 
implies that
the spaces $\mathcal{M}_{c}$ are in fact orbit spaces of certain $K(\pi ,1)$%
-spaces under an action of \textit{Diff}$(S^{1})$. This result is proved
explicitly in \cite{hatcher2}.

It is a result of Gluck and Pan \cite{gluckandpan}
that the set of knots in any fixed
isotopy class with non-vanishing curvature and a prescribed self-linking
number is path-connected. In a future paper, we will show that the inclusion
of the set of all such knots into the set of all knots in a particular
isotopy class is also a weak homotopy equivalence.

\end{document}